\input amstex
\input amsppt.sty
\documentstyle{amsppt}
\topmatter 
\title On $\pi - \pi$ theorem for manifold pairs with boundaries
\endtitle
\author Matija Cencelj -- Yuri V. Muranov -- Du\v san Repov\v s 
\endauthor
\subjclass  Primary  57R67, 57Q10
Secondary  57R10, 55U35
\endsubjclass
\thanks{Second author partially supported by the Russian 
Foundation for Fundamental 
Research  Grant  05--01--00993, the first and third authors 
partially supported by the Ministry of Higher Education, Science and 
Technology of the Republic of Slovenia research program P1--292--0101-04.}
\endthanks
\keywords
Surgery  obstruction groups, normal map, homotopy triangulation, 
surgery on manifold pairs, 
splitting obstruction groups, $\pi$-$\pi$ theorem
\endkeywords
\abstract Surgery obstruction of a normal 
map to a 
simple Poincare pair $(X,Y)$ lies in the relative surgery obstruction 
group $L_*(\pi_1(Y)\to\pi_1(X))$.    
A well known result of Wall, 
the so called  $\pi$-$\pi$ theorem,   states that in higher dimensions 
a normal map of a manifold with boundary to a simple Poincare pair
 with $\pi_1(X)\cong\pi_1(Y)$ is normally bordant
to a simple homotopy equivalence of pairs. 
In order to
study normal maps to a manifold with  a 
submanifold, 
Wall introduced surgery obstruction group for manifold pairs
$LP_*$  
and splitting obstruction groups $LS_*$.           
In 
the present paper we formulate and prove
for manifold pairs with boundaries 
the results which are 
similar to the $\pi$-$\pi$ theorem. 
We give  
direct geometric proofs, which are based on the original 
statements of Wall's results  and 
apply obtained results to investigate
surgery on filtered 
manifolds.  
\endabstract
\endtopmatter

\document

\bigskip
        
\subhead 1. Introduction 
\endsubhead

The surgery obstruction groups 
$L_{*}(\pi)$ were introduced by Wall in his fundamental paper \cite{7} . 
Let $(f, b) : M\to X$ be a normal map from a closed manifold $M^n$  
to a  simple Poincar\'e complex $X$ of formal dimension $n$ where 
 $b : \nu_M \to \xi$ is a bundle
map which covers $f : M \to X$. Then an  obstruction
$\theta(b, f)\in L_n(\pi_1(X))$ for the existence of a simple homotopy
equivalence in the class of normal bordism  of the
map $(f,b)$ is defined. 

Indeed,
Wall 
defines $L_*$-groups and  describes surgery 
theory for the case of 
manifold
$n$-ads.
For example,  an obstruction to surgery on the manifold 
pair  
lies in the relative surgery obstruction group 
$L_*(\pi_1(Y)\to\pi_1(X))$.
Hence,  if the map 
$\pi_1(Y)\to\pi_1(X)$ is an isomorphism, then   
in high dimensions 
a normal map from a manifold with boundary to
 the  simple Poincare pair
 is normally bordant
to a simple homotopy equivalence of pairs. 
In \cite{7, \S 4} Wall gives direct
geometric
proof of this 
"important special case". 

The surgery obstruction groups
and natural maps
do not depend on
the category of manifolds (see \cite{6} and \cite{7}). 
In the present  paper we shall work in the category of
topological manifolds. All results of 
this  paper are 
transferred to $Diff$  or $PL$@-manifolds.

Let $(X^n, Y^{n-q}, \xi)$ be a codimension $q$ manifold pair
\cite{6, \S 7.2}.
Let
$$
F=\left(\matrix {\pi}_{1}(S(\xi)) & \longrightarrow &
{\pi}_{1}(X\setminus Y) \cr
\downarrow & \  & \downarrow \cr
{\pi}_{1}(Y) & \longrightarrow & {\pi}_{1}(X)
\endmatrix \right),
\tag 1
$$
be the pushout square  of fundamental groups with 
orientations 
where  $S(\xi)$ is 
the  boundary of the tubular neighborhood
of $Y$ in $X$. 
 
The obstruction groups $LP_{*}(F)$ to surgery on the manifold 
pair (to obtain $s$@-triangulation of the pair) 
are defined (see \cite{7, \S 17E} and \cite{6, \S 7.2}). 
A simple homotopy equivalence $f: M\to X$ splits along a 
submanifold   $Y$ if it is homotopy equivalent to a map $g$ which is 
an $s$@-triangulation of the manifold pair $(X, Y, \xi)$. 
The splitting obstruction groups 
$LS_{*}(F)$ 
are defined in \cite{6} and \cite{7}.

Let 
$$
(Y,\partial Y)\subset (X,\partial X)
\tag 2
$$
be a codimension $q$ manifold pair with boundaries (see  \cite{6, \S 7}). 
A normal map 
$$
(f,\partial f): (M,\partial M)\to (X, \partial X)
\tag 3
$$
is an $s$@-triangulation of the manifold pair $(X, Y)$ 
with boundaries $(\partial X, \partial Y)$  
if the maps 
$
f: M\to X
$
and 
$
\partial f: \partial M\to \partial X
$
are $s$@-triangulations of the pairs $(X,Y)$ and 
 $(\partial X, \partial Y)$,  respectively.  

In $rel_{\partial}$  case we 
 consider  normal maps   (3)
which are already 
split
on the boundary. 
The surgery theory 
for this case was developed
in \cite{6} and \cite{7}. 
In the $rel_{\partial}$@-case the obstructions to surgery on manifold 
pairs relatively the boundary lie in the group $LP_*(F)$, and similarly 
for  $rel_{\partial}$ splitting obstruction (see \cite{6, pages 584--587}).

In the present paper  we consider  surgery 
on manifold pairs with boundaries without
fixing of maps on the 
boundary. To study surgery on filtered manifolds (see \cite{1}, \cite{5}, 
and \cite{8}) we need to know 
the 
geometric
properties of surgery on manifold 
pairs with boundaries. This is the  special 
case of splitting theory 
for manifolds $n$@-ads mentioned by Wall 
\cite{7, page 136}. We give 
the exact statement and proof  of $\pi - \pi$@-theorem  for  
various  maps to  a manifold pair with 
boundary.  Then we
apply 
obtained  results to investigation of surgery on  filtered  manifolds.

A manifold pair (2)
with 
boundary
defines a  pair of closed manifolds
$\partial Y\subset \partial X$ with 
 a pushout square   
$$
\Psi=\left(\matrix {\pi}_{1}(S(\partial\xi) & \longrightarrow &
{\pi}_{1}(\partial X\setminus \partial Y) \cr
\downarrow & \  & \downarrow \cr
{\pi}_{1}(\partial Y)) & \longrightarrow & {\pi}_{1}(\partial X)
\endmatrix \right)
\tag 4
$$
of fundamental groups for the splitting problem. Here $\partial\xi$ 
is the restriction of $\xi$ on the boundary $\partial Y$.

A natural inclusion $\delta: \partial X \to X$ induces 
a map of $\Delta:\Psi \to F$ of squares of fundamental groups. 

Now we can formulate the main result of the paper.

\proclaim{Theorem}  Let 
$
(Y^{n-q},\partial Y)\subset (X^n,\partial X)
$
be a codimension $q$ manifold pair with
boundary
with $n-q\geq 6$. 
Let the map $\Delta$ be an isomorphism 
of the squares of fundamental groups. Under these conditions  
we have the following results.

A) Any  normal  map (3)
is  normally bordant to an $s$-triangulation of  
$
(X, \partial X)
$
which is split along $(Y, \partial Y)\subset (X, \partial X)$.

B) Any  simple homotopy equivalence (3) of pairs 
is concordant  to    
a simple homotopy equivalence to   
$
(X, \partial X)
$
which is split along $(Y, \partial Y)\subset (X, \partial X)$.

C) Let a normal map (3)
define
a simple homotopy equivalence of pairs  
$$
f|_{(N, \partial N)}: (N, \partial N)\to (Y, \partial Y)
$$
where $N=f^{-1}(Y),\ \partial N=f^{-1}(\partial Y)$ 
are transversal preimages. Then $(f,\partial f)$
is normally bordant  
to an $s$@-triangulation of $(X, \partial X)$ which is split 
 along $(Y, \partial Y)\subset (X, \partial X)$.
Moreover, there 
exists
a transversal to $Y\times I$ bordism 
$$
F=(h;g, f_0, f_1):
(W; V, M, M^{\prime})\to (X\times I; \partial X\times I,  X\times\{0\},
 X\times\{1\})
$$
for which the restriction 
$
h_{h^{-1}(Y\times I)}
$
is 
$$
(f|_{(N, \partial N)})\times \operatorname{Id}: 
(N, \partial N)\times I\to (Y, \partial Y)\times I
$$ 
\endproclaim 

In section 2 we give necessary preliminary material. In section 3 we 
prove
the theorem and  apply our results to 
surgery on filtered manifolds. 
\smallskip

\subhead 2. Preliminaries
\endsubhead
\bigskip

We shall consider a case of topological manifolds  and follow
notations from \cite{6, \S 7.2}. 
Let $(X, Y, \xi)$ be a codimension $q$   manifold 
pair  in the sense of Ranicki (see \cite{6, \S 7.2}) i. e.  
a locally flat closed submanifold  $Y$ is given with a normal fibration 
 $$ 
  \xi=\xi_{Y\subset X}: Y \to \widetilde{BTOP}(q) 
 $$ 
with the associated $(D^q, S^{q-1})$ 
fibration 
$$
(D^q, S^{q-1})\to (E(\xi), S(\xi))\to Y
$$ 
and we have a decomposition of the closed manifold  
$$ 
 X =E(\xi)\cup_{S(\xi)}\overline{X\setminus E(\xi)}. 
 \tag 5
$$ 
A topological normal map  \cite{6, \S 7.2}
$$ 
((f,b), (g,c)):(M,N)\to (X,Y) 
$$
to the manifold pair $(X, Y, \xi)$) is 
represented by a normal map 
$(f,b)$ to the manifold $X$ 
which is transversal to $Y$ with $N=f^{-1}(Y)$, and 
$(M,N)$ is a topological manifold pair with a normal 
fibration 
$$ 
\nu:N\overset{f|_N}\to{\to} Y \overset{\xi}\to{\to} 
\widetilde{BTOP}(q). 
$$ 
Additionaly the following conditions are satisfied:

(i) the restriction 
$$ 
(f,b)|_N =(g,c) :N\to Y 
$$ 
is a normal map;

(ii) the restriction 
$$ 
(f,b)|_P =(h,d) :(P, S(\nu))\to (Z,S(\xi)) 
$$ 
is a normal map to the pair $(Z,S(\xi))$, 
where 
$$ 
P=\overline{M\setminus E(\nu)}, \ \ Z=\overline{X\setminus 
E(\xi)}; 
$$ 

(iii) the restriction 
$$ 
(h,d)|_{S(\nu)}: S(\nu)\to S(\xi) 
$$ 
coincides with the induced map 
$$ 
(g,c)^!: S(\nu)\to S(\xi), 
$$ 
and $(f,b)=(g,c)^!\cup(h,d)$. 

The normal maps to $(X,Y,\xi)$ are called $t$-triangulations 
of the manifold pair $(X,Y)$.
Note
that the set of concordance classes of $t$-triangulations 
of the pair $(X, Y, \xi)$ coincides with the set of $t$-triangulations 
of the manifold $X$ \cite{6, Proposition 7.2.3}. 

 An $s$@-triangulation of a manifold pair $(X,Y, \xi)$  in 
topological category 
 \cite{6, p. 571}
is a $t$@-triangulation 
 of this pair for which the maps 
 $$ 
 f:M\to X,\ g:N\to Y, \ 
 \text{and} \ 
 (P, S(\nu))\to (Z, S(\xi)) 
\tag 6
 $$ 
are simple homotopy equivalences ($s$@-triangulations).

For a 
codimension $q$ manifold pair with boundaries (2) 
 \cite{6, p. 585}
we  have a normal fibration 
$(\xi, \partial \xi)$ over the pair $(Y, \partial Y)$
and, similarly to (5),  a decomposition  
$$
(X, \partial X)= (E(\xi)\cup_{S(\xi)}Z, 
E(\partial\xi)\cup_{S(\partial\xi)}\partial_+Z)
\tag 7
$$
where $(Z; \partial_+Z, S(\xi); S(\partial\xi))$ is a manifold
triad. Note here that $\partial_+Z =\overline{\partial X\setminus 
E(\partial\xi)}$.

A topological 
normal map (4) 
of manifold pairs with boundaries provides 
a normal fibration 
$(\nu, \partial \nu)$ over the pair $(N, \partial N)$
where $(N, \partial N)=(f^{-1}(Y),(\partial f)^{-1}(\partial Y))$ 
\cite{6, p. 570}. We have the 
following decomposition 
$$
(M, \partial M)= (E(\nu)\cup_{S(\nu)}P, 
E(\partial\nu)\cup_{S(\partial\nu)}\partial_+P)
\tag 8
$$
where $(P; \partial_+P, S(\nu); S(\partial\nu))$ is a manifold
triad.  
Recall, that two $s$@-triangulations of the pair $(X, \partial X)$
$$
(f_i, \partial f_i): (M_i, \partial M_i)\to (X, \partial X), \ i=0,1
$$
are concordant (see \cite{7, \S 10} and \cite{6, \S 7.1})
if there exists a simple homotopy equivalence of $4$@-ads
$$
(h;g, f_0, f_1):
(W; V, M_0, M_1)\to (X\times I; \partial X\times I,  X\times\{0\},
 X\times\{1\})
\tag 9
$$
with
$$
\partial V=\partial M_0\cup \partial M_1.
$$

\smallskip

\subhead 3. Proof  and Corollary
\endsubhead 
\bigskip

 Consider the case A) of the theorem. 
 A restriction of the map 
$$
(f,\partial f): (M,\partial M)\to (X, \partial X)
$$
gives   
a normal map of pairs  
$$
\phi=f|_{S(\nu)}: (S(\nu), S(\partial \nu))\to(S(\xi), S(\partial \xi)).
\tag 10
$$ 
The isomorphism $\Delta$ provides an 
isomorphism 
$\pi_1(S(\partial \xi))\to \pi_1(S(\xi))$.
Hence the normal map $\phi=f|_{S(\nu)}$ (10)
of the 
pairs
satisfies the conditions of 
$\pi - \pi $ theorem of Wall \cite{7, \S 4} and 
it is normally 
bordant to a simple homotopy equivalence of pairs. By 
\cite{3, page 45} we can extend this bordism to obtain a normal bordism
with the bottom map $(f, \partial f)$ and with the top
a normal map 
$$
(f^{\prime}, \partial f^{\prime}) : 
(M^{\prime},\partial M^{\prime})\to (X, \partial X)
$$
with the properties similar to the map $f$  and for which the 
restriction 
$$
\phi^{\prime}=f^{\prime}|_{S(\nu^{\prime})}: 
(S(\nu^{\prime}), S(\partial \nu^{\prime}))
\to (S(\xi), S(\partial \xi)).
$$
is a simple homotopy equivalence of pairs. 
To avoid complicated notations  we can suppose
that the normal map (4) has the restriction 
$\phi$ (10) which is a  simple homotopy equivalence 
of pairs. 

Now the  restriction of the map $(f, \partial f)$ gives   
a normal map of triads   
$$
\psi = f|_{E(\nu)}: (E(\nu); E(\partial \nu), S(\nu); S(\partial\nu)) \to
(E(\xi); E(\partial \xi), S(\xi); S(\partial \xi)).
$$ 
The restriction 
$$
\psi|_{S(\nu)}=\phi :(S(\nu), S(\partial\nu)) \to (S(\xi), S(\partial \xi))
$$
is the simple homotopy equivalence of pairs, and 
the isomorphism $\Delta $ provides the
isomorphism   
$\pi_1(E(\partial \xi))\to \pi_1(E(\xi))$.
Hence the normal map $\psi $ satisfies the conditions of 
$\pi - \pi $ theorem for the triad  \cite{7, Theorem 3.3}
and  it is  normally 
bordant to a simple homotopy equivalence of triads by bordism
$\Phi: V\to E(\xi)\times I$  
relatively $(S(\nu), S(\partial \nu))$
with the bottom map 
$\psi : V_0=E(\nu)\to E(\xi)$. 
The restriction of 
$(\Phi,V)$ to the $S(\nu)$ is a trivial bordism
$$
\phi_t: Q=S(\nu)\times I\to S(\xi)\times I, \ \phi_t(x, t)= (\phi(x), t), \
  t\in  I=[0,1].
$$    

We can attach the bordism  $V$ to the manifold $M$ 
identifying $V_0$ with $E(\nu)$ to obtain 
a space $\Lambda= V\cup_{E(\nu)}M$.     
Since the restriction of $\Phi$ to the
bottom $V_0=E(\nu)$ coincides with the map $\psi=f|_{E(\nu)}$
we obtain the map  
$$
F=\Phi\cup_{\psi}f:\Lambda =V\cup_{E(\nu)}M\to (E(\xi)\times I)
\cup_{E(\xi)\times\{0\}} X.
$$

In a similar way,  the restriction of the map $(f, \partial f)$ gives   
a normal map of triads   
$$
\alpha = f|_{P}: (P; \partial_+P, S(\nu); S(\partial\nu)) \to
(Z; \partial_+Z, S(\xi); S(\partial \xi)), 
$$ 
for which the restriction  
$$
\alpha|_{S(\nu)}=\phi:
(S(\nu), S(\partial\nu)) \to (S(\xi), S(\partial \xi))
$$
is a simple homotopy equivalence of pairs.
The isomorphism $\Delta$ provides an 
isomorphism   
$\pi_1(\partial_+Z)\to \pi_1(Z)$.
Hence the normal map $\alpha $ is  normally 
bordant to a simple homotopy equivalence of triads by bordism  
$G:W\to Z\times I$  
relatively  $(S(\nu), S(\partial \nu))$ with the bottom map 
$$
G|_{W_0}=\alpha=f|_P :W_0=P \to Z.
$$ 

We can attach $W$ to 
the space $\Lambda$  to obtain 
a space $\Lambda\cup_{Q\cup P} W$.
The restriction $G|_{W_0}$ coinsides with the map $\alpha = f|_{P}$
and bordism maps $F$ and $G$ coincides on $Q$.  
Thus we obtain a bordism 
$$
\Omega: \Lambda\cup_{Q\cup P} W\to X\times I
\tag 11
$$
where the map $\Omega$  extends the maps $F$ and $G$.
By our construction,  on the top of 
the bordism (11) we obtain the map 
 $$
(f^{\prime}, \partial f^{\prime}) : 
(M^{\prime},\partial M^{\prime})\to 
(X, \partial X)
$$
for which the restrictions give 
the  simple homotopy equivalences 
of triads 
$$
(E(\nu^{\prime}); E(\partial \nu^{\prime} ), 
S(\nu ^{\prime}); S(\partial\nu ^{\prime})
\to 
(E(\xi); E(\partial \xi), S(\xi); S(\partial \xi))
$$
and 
$$
(P^{\prime}; \partial_+P^{\prime}, 
S(\nu ^{\prime}); S(\partial\nu^{\prime}))
\to 
(Z; \partial_+Z, S(\xi); S(\partial \xi)), 
$$
where 
$$
(M^{\prime}, \partial M^{\prime})=
(E(\nu^{\prime})\cup_{S(\nu^{\prime})}P^{\prime}, 
E(\partial \nu^{\prime})\cup_{S(\partial\nu^{\prime})}P_+^{\prime}).
$$ 

To finish the proof of A) we must verify  only that the constituent 
maps $f^{\prime}:M^{\prime}\to X$ and 
$\partial f^{\prime}: \partial M^{\prime}\to \partial X $ 
are simple homotopy equivalences.

The space  $M^{\prime}$ is union of two parts
$$
P^{\prime}\cup_{S(\nu^{\prime})}E(\nu^{\prime})
$$ 
with the
intersection  $S(\nu ^{\prime})$. The  restrictions of $f$ on these
two parts and on the intersection are simple homotopy equivalences. 
Hence by \cite{2, Theorem 23.1} the map 
$f^{\prime}:M^{\prime}\to X$ is simple homotopy equivalence.
For the map 
$\partial f^{\prime}$ the situation is similar
since  
$$
\partial M^{\prime}=\partial_+P^{\prime}
\cup_{S^{\prime}(\partial\nu)}E^{\prime}(\partial\nu).
$$
The case A) is proved.
 
Now consider the case B). 
The map $(f,\partial f)$ (4) is a simple homotopy equivalence 
of pairs. Since the map $\Delta $ is an isomorphism then by A) the map 
 $(f,\partial f)$  is normally bordant to a map 
$$
(f^{\prime}, \partial f^{\prime}) : 
(M^{\prime},\partial M^{\prime})\to 
(X, \partial X)
$$
which is split along $(Y, \partial Y)\subset (X, \partial X)$.
Thus we have a normal bordism 
$$
\Phi: W\to X\times I, \  \partial W= W_0\cup W_1\cup V
$$
where 
$$
W_0=M, W_1= M^{\prime}, V= \Phi^{-1}(\partial X\times I).
$$

The bordism $(W, \Phi)$ gives a normal 
map of  manifold triads 
$$
(W; V, W_0\cup W_1; \partial W_0\cup \partial W_1)
\to
(X\times I; \partial X\times I, X\times \{0,1\}; \partial X\times\{0,1\})
\tag 12
$$
which we shall denote by $\Phi$, too. 

The restriction of the normal map $\Phi$ on  
$W_0\cup W_1$ 
is the simple homotopy equivalence of pairs 
$$
(W_0\cup W_1, \partial W_0\cup \partial W_1) \to 
(X\times \{0,1\}; \partial X\times\{0,1\}).
$$
The fundamental group of the triad 
$$
(X\times I; \partial X\times I, X\times \{0,1\}; \partial X\times\{0,1\})
$$
is equal to 
$$
F_I=\left(\matrix {\pi}_{1}(\partial X\times\{0,1\})) & \longrightarrow &
{\pi}_{1}(\partial X\times I) \cr
\downarrow & \  & \downarrow \cr
{\pi}_{1}(X\times\{0,1\})) & \longrightarrow & {\pi}_{1}(X\times  I)
\endmatrix \right).
\tag 13
$$
The isomorphism  $\Delta$ provides the 
 isomorphism $\pi_1(\partial X)\to
\pi_1(X)$ and,  hence,  the  vertical maps 
in the  square (13) are isomorphisms of grouppoids. 
Hence the normal map of triads (12) 
satisfies the conditions 
of $\pi-\pi$ theorem for triads relatively 
$$ 
\partial_2(X\times I)=
\left((\partial X\times\{0,1\})\subset (X\times \{0,1\})\right).
$$
Thus the map $\Phi$ (12) of triads 
is normally bordant relatively  $\partial_2(X\times I)$
to a simple homotopy equivalence of triads 
$$
\Phi^{\prime}:(W^{\prime}; V^{\prime}, W_0\cup W_1; \partial W_0\cup \partial W_1)
\to
(X\times I; \partial X\times I, X\times \{0,1\}; \partial X\times\{0,1\})
$$
where 
$$
\Phi^{\prime}|_{W_0}=f, \ \Phi^{\prime}|_{W_1}=f^{\prime}, \
\Phi^{\prime}|_{\partial W_0}=\partial f, \ 
\Phi^{\prime}|_{\partial W_1}=\partial f^{\prime}.
$$

By our construction the map $(f^{\prime}, \partial f^{\prime})$  is splitted
along the pair  $(Y, \partial Y)\subset (X, \partial X)$ and by (9) the map 
$\Phi^{\prime}$ gives the concordance between $(f, \partial f)$ and 
$(f^{\prime}, \partial f^{\prime})$.
The part B) of the theorem is proved.  

In the case C) the map $f$ is transversal to $(Y, \partial Y)$ and its
restriction is  a simple homotopy equivalence 
$(N, \partial N)\to (Y, \partial Y)$. Hence $f|_N$ induces a simple homotopy 
equivalence of tubular neighborhoods
with boundaries (see \cite{7, page 8} and \cite{6, page 579}).
We obtain a simple homotopy equivalence of 
triads 
$$
\psi:  (E(\nu); E(\partial \nu), S(\nu); S(\partial\nu)) \to
(E(\xi); E(\partial \xi), S(\xi); S(\partial \xi))
$$ 
for which  the restriction  
$$
\psi|_{(N, \partial N)}=f|_{(N, \partial N)}
$$
is the simple homotopy equivalence of pairs. Now the result follows 
from consideration 
of the map $\alpha $ of triads from A) by the same arguments. 
The theorem is proved.
\qed 
\smallskip

Now  we 
apply obtained results 
to  surgery on filtered manifolds 
(see \cite{1}, \cite{4}, \cite{5},
and \cite{8}). At first we recall necessary definitins.

Let  $Z^{n-q-q^{\prime}}\subset Y^{n-q} \subset X^n$ be a triple 
of closed topological manifolds (see \cite{4}, \cite{5}, and \cite{6}).
We have the following topological normal bundles:
  $\xi$  for the submanifold $Y$ in $X$,
  $\eta$ for the submanifold  $Z$ in  $Y$, and
  $\nu$ for the submanifold $Z$ in $X$.
Let $(E(\xi),S(\xi))$, $(E(\eta), S(\eta))$, and $(E(\nu), S(\nu))$, 
respectively,
be the spaces with boundaries  
of associated $(D^*, S^{*-1})$ fibrations. 
We identify the space $E(\nu)$ 
with the space  $E^{\prime}(\xi)$
of the restriction
 $\xi|_{E(\eta)}$ in such a way 
 that 
$$
S(\nu)=E^{\prime\prime}(\xi)\cup S^{\prime}(\xi)
\tag 14
$$
where $E^{\prime\prime}(\xi)$ is the space of the restriction 
$\xi|_{S(\eta)}$ and $S^{\prime}(\xi)$ is the restriction 
of $S(\xi)$ on  $E(\eta)$  (see \cite{1} and \cite{5}).

Let 
$$
(X_k, \partial X_k)\subset (X_{k-1}, \partial X_{k-1}) \subset \cdots 
\subset (X_0, \partial X_0) =(X, \partial X)
\tag 15
$$
be a filtration  of a compact manifold 
$(X, \partial X)$
by manofolds with boundaries  (see \cite{1}, \cite{5}, and \cite{8}). 
From now we shall assume that the dimension 
$\operatorname{dim} X_k\geq 6$.
  
The filtration (15) defines the filtration 
$$
\partial X_k\subset  \partial X_{k-1} \subset \cdots 
\subset \partial X_0 =\partial X
\tag 16
$$
of $\partial X$ by closed manifolds.  
Recall that any triple of manifolds from (15) and (16) satisfy 
properties that are similiar to (14) on the corresponding 
normal bundles. Additionaly,  for every pair of manifolds 
with boundaries from (15) we have a decomposition 
that is similiar to (7).  The filtration (15) 
defines a stratified manifold with boundary 
$(\Cal X, \partial \Cal X)$ (see \cite{1}, \cite{5}, and \cite{8}). 

Any topological normal map (3)  
defines the topological normal map to the filtration  (15) 
(see \cite{5} and \cite{8}). 
Let $M_i$ be the transversal preimage of the submanifold 
$X_i$.
A topological normal map (3)
is an $s$@-triangulation of the filtration (15)
if constituent normal maps 
$$
f|_{(M_j, M_i)}:(M_j, M_i)\to (X_j, X_i), \ 0\leq j\le i\leq k 
$$ 
are $s$@-triangulations of the manifold pairs with boundaries   $(X_j, X_i)$. 
The stratified Browder-Quinn surgery obstruction
groups $L_n^{BQ}(\Cal X, \partial \Cal X)$ are defined
\cite{1}.

For $(1\leq i\leq k)$,  let $F_i$   
be the square in the splitting problem for the 
manifol pair $X_i\subset X_{i-1}$,  and $\Psi_i$ be the similiar 
square for the closed manifold pair 
$\partial X_i\subset \partial X_{i-1}$. The natural inclusions 
of boundaries induce the maps 
$$
\Delta_i:\Psi_i\to F_i
$$
 for $1\leq i\leq k$.  
\smallskip

\proclaim{Corollary}  Let all the maps $\Delta_i$ $(1\leq i\leq k)$ 
be
isomorphisms.  Then every normal map to the filtered space $\Cal X$ 
is normally bordant to an $s$@-triangulation of $\Cal X$ and, 
hence,  the group  $L_n^{BQ}(\Cal X, \partial \Cal X)$ is trivial.
  \endproclaim

\demo{Proof} Denote by $f_{k-1}$  the restriction of the normal map $f$ on
the manifold  $M_{k-1}$.  By item A)  of 
the theorem the map 
$f_{k-1}$ is normally bordant to an $s$@-triangulation $g_{k-1}$ of 
the pair $(X_{k-1}, X_k)$. By \cite{3} we 
can extend this bordism to 
obtain a bordism $F: W\to X\times I$  
with a top  normal map 
$
(g, \partial g)
$
to $(X, \partial X)$
whose restriction on $M_{k-1}$ is the $s$@-triangulation $g_{k-1}$.

Consider the resriction 
$(g, \partial g)|_{M_{k-2}}$ for which the the restriction 
on $M_{k-1}$ is the $s$@-triangulation $g_{k-1}$.
By item C) of 
the theorem the map
$(g, \partial g)|_{M_{k-2}}$ is normally bordant to a map 
$g_{k-2}$ which is an $s$@-triangulation of the pair 
$(X_{k-2}, X_{k-1})$ and 
$$
g_{k-2}\left|_{M_{k-1}}\right.=g_{k-1}.
$$ 

Now, 
by applying item C) $(k-1)$ time, we obtain a map
$$
(g_0, \partial g_0):(M^{\prime}, \partial M^{\prime})\to (X, \partial X)
$$
with the following properties. For any $0\leq j\leq k-1$ 
the restriction of  $(g_0, \partial g_0)$ on the transversal preimage 
of $(X_j, X_{j+1})$ is the $s$@-triangulatin $g_j$ of the manifold 
pair with boundaries  $(X_j, X_{j+1})$.  
If a normal map $f:M\to X$ is an $s$@-triangulation 
of the subfiltration  
$$
X_{k-1}\subset \cdots \subset X_{1}\subset X
\tag 17
$$
and the  restriction $f\left|_{M_{k}}\right.$ is an $s$-triangulation 
of the pair $(X_{k-1}, X_{k})$, 
 then  by \cite{5, Proposition 2.5}, the map $f$ is an 
 $s$@-triangulation of the filtartion (15) $\Cal X$.
Now we can 
apply this result $k-1$ times  starting from 
the subfiltration $X_1\subset X_0$ of the filtration  
$X_2\subset X_1\subset X_0$ untill  the   subfiltration 
(17) of the filtartion (15). Corollary is proved. \qed 
\enddemo

\vfill\eject

Authors' addresses:

\bigskip

\noindent
Matija Cencelj :

Institute for Mathematics, Physisc and Mechanics, University of 
Ljubljana, Jadranska 19, Ljubljana, Slovenia 
\noindent
email: matija.cencelj\@fmf.uni-lj.si

\bigskip

\noindent
Yuri V. Muranov:  Department of Information Science and Management, 
Institute of Modern Knowledge, 
ulica Gor'kogo 42,
210004 Vitebsk,
Belarus;

\noindent
email: ymuranov\@imk.edu.by

\bigskip

\noindent
Du\v san Repov\v s:

Institute for Mathematics, Physisc and Mechanics, University of 
Ljubljana, Jadranska 19, Ljubljana, Slovenia 
\noindent
email: dusan.repovs\@fmf.uni-lj.si
\enddocument

\bye